\documentclass[11pt]{amsart}
\usepackage{amsmath}
\usepackage{amsfonts}
\usepackage{graphicx,color, xcolor}
\usepackage{amssymb}
\usepackage{amsthm}
\usepackage{comment}
\usepackage[left =1 in , right = 1 in, bottom = 1 in, top = 1 in]{geometry}
\usepackage[all]{xy}

\makeatletter
\newtheorem*{rep@theorem}{\rep@title}
\newcommand{\newreptheorem}[2]{
\newenvironment{rep#1}[1]{
\def\rep@title{#2 \ref{##1}}
\begin{rep@theorem}}
{\end{rep@theorem}}}
\makeatother

\theoremstyle{plain}
\newtheorem{thm}{Theorem}[section]
\newtheorem{cor}[thm]{Corollary}
\newtheorem{lem}[thm]{Lemma}
\newtheorem{prop}[thm]{Proposition}

\newreptheorem{theorem}{Theorem}
\newreptheorem{cor}{Corollary}

\theoremstyle{definition}
\newtheorem{de}[thm]{Definition}

\theoremstyle{remark}
\newtheorem{rem}[thm]{Remark}
\newtheorem{example}[thm]{Example}

\newcommand{\Spec}{\operatorname{Spec}}

\newcommand{\Hom}{\operatorname{Hom}}
\newcommand{\Rat}{\operatorname{Rat}}

\newcommand{\SL}{\operatorname{SL}}
\newcommand{\PGL}{\operatorname{PGL}}
\newcommand{\Stab}{\operatorname{Stab}}

\newcommand{\NS}{\operatorname{NS}}
\newcommand{\Pic}{\operatorname{Pic}}

\newcommand{\Aut}{\operatorname{Aut}}

\newcommand{\Sym}{\operatorname{Sym}}

\newcommand{\codim}{\operatorname{codim}}
\newcommand{\M}{\operatorname{M}}

\newcommand{\mc}[1]{\mathcal{#1}}
\newcommand{\mb}[1]{\mathbb{#1}}

\newcommand{\ti}[1]{\textit{#1}}
\newcommand{\tr}[1]{\textrm{#1}}

\newcommand{\bd}{\begin{de}}
\newcommand{\ed}{\end{de}}
\newcommand{\bl}{\begin{lem}}
\newcommand{\el}{\end{lem}}
\newcommand{\bp}{\begin{prop}}
\newcommand{\ep}{\end{prop}}
\newcommand{\bt}{\begin{thm}}
\newcommand{\et}{\end{thm}}
\newcommand{\bc}{\begin{cor}}
\newcommand{\ec}{\end{cor}}
\newcommand{\bpf}{\begin{proof}}
\newcommand{\epf}{\end{proof}}

\newcommand{\beq}{\begin{equation}}
\newcommand{\eeq}{\end{equation}}
\newcommand{\beqs}{\begin{equation*}}
\newcommand{\eeqs}{\end{equation*}}
\newcommand{\ben}{\begin{enumerate}}
\newcommand{\een}{\end{enumerate}}
\newcommand{\bit}{\begin{itemize}}
\newcommand{\eit}{\end{itemize}}

\title{Isotriviality and the Space of Morphisms on Projective Varieties}

\author{Anupam Bhatnagar}
\address{Department of Mathematics, Borough of Manhattan Community College, The City University of New York; 199 Chambers Street, New York, NY 10007 U.S.A.}
\email{abhatnagar@bmcc.cuny.edu}

\author{Alon Levy}
\address{Department of Mathematics, The University of British Columbia; 1984 Mathematics Road
Vancouver, B.C. Canada V6T 1Z2}
\email{levy@math.ubc.ca}

\date{\today}

\begin{document}

\begin{abstract}
Let $K=k(C)$ be the function field of a smooth projective curve $C$ over an infinite
field $k$, let $X$ be a projective variety over $k$. We prove two results. First, we
show with some conditions that a $K$-morphism $\phi: X_K \to X_K$ of degree at least two is isotrivial
if and only if $\phi$ has potential good reduction at all places $v$ of $K$. Second,
let $(X,\phi), (Y,\psi)$ be dynamical systems where $X,Y$ are defined over $k$
and $g:X_{K} \to Y_{K}$ a dominant $K$-morphism, such that $g \circ \phi = \psi \circ g$.
We show under certain conditions that if $\phi$ is defined over $k$, then $\psi$ is
defined over $k$.
\end{abstract}

\maketitle

\section{Introduction}
Let $C$ be a smooth projective curve over an infinite field $k$ and let $K=k(C)$ denote the
function field of $C$. In \cite{PST} the authors prove a criterion for isotriviality of endomorphisms
on $\mb{P}^n_K$ in terms of good reduction at every place of $K$. In this paper we generalize this
criterion to projective $K$-varieties. To state our result we need a few definitions.

\bd \label{isotrivial}
A projective $K$-variety $X$ is \ti{trivial} if it is defined over the constant
field $k$ and $X$ is \ti{isotrivial} if there exists a finite extension $K'/K$ such that
$X_{K'}:= X \times_{\Spec(K)} \Spec(K')$ is trivial. Now let $X$ be a projective $k$-variety.
A morphism
$\phi: X_K \to X_K$ is \ti{trivial} if it is defined over the constant field $k$ and
it is \ti{isotrivial} if there exists a finite extension $K'/K$ such that the induced
morphism $\phi': X_{K'} \to X_{K'}$ is defined over $k'$, the algebraic closure
of $k$ in $K'$.
\ed

Let $\mc{X}$ be a projective model of $X_K$ over $C$ and for some non-empty open
subset $U$ of $C$, let $\mc{X}_U: =\mc{X} \times_C U$. Let $M_K$ denote the places
of $K$ and identify the places of $K$ with the closed points of $C$.

\bd \label{pgr}
The morphism $\phi :X_K \to X_K$ has \ti{good reduction} at $v$ if there exists an
open subset $U\subset C$ containing $v$ such that $\phi$ extends to a morphism
on $\mc{X}_U$. We say $\phi$ has \ti{potential good reduction} at $v$ if there exists a
finite extension $K'$ of $K$ and a place $v'$ of $K'$ over $v$ such that $\phi'$ has
good reduction at $v'$.
\ed

\bd Let $X$ be a projective variety and $\phi: X\to X$ a morphism. The dynamical system 
$(X,\phi)$ is said to be \ti{polarized} if for some ample line bundle $\mc{L}$ on $X$, we 
have $\phi^*\mc{L} \cong \mc{L}^{\otimes q}$ with $q>1$.
\ed

In \S 2 we prove our first result, which states:

\begin{repcor}{isotriv}
Let $K=k(C)$ be the function field of a smooth projective curve $C$ defined over
an infinite field $k$, $X$ a connected projective $k$-variety, $\mc{L}$ an ample line bundle on $X$,
$\phi: X_K \to X_K$ a polarized $K$-morphism. Assume that
the subgroup of automorphisms preserving the line bundle class
$$\Aut(X, \mc{L}):= \{ \tau \in \Aut(X)\; |\; \tau^*\mc{L} \cong \mc{L} \} $$
is reductive. Then $\phi$ is isotrivial if and only if $\phi$ has potential good
reduction at all places $v$ of $K$.
\end{repcor}

The proof of the theorem uses geometric invariant theory(GIT) and our
approach is similar to \cite{PST} (Thm. 1). The condition that $\Aut(X,\mc{L})$
is reductive is nontrivial; there exist varieties for which $\Aut(X,\mc{L})$ is non-reductive
that nonetheless admit polarized endomorphisms. The example with principally polarized
abelian varieties in \S1 of \cite{PST} shows that one cannot have a result similar to
Corollary \ref{isotriv} for an arbitrary $K$-variety.

In \S 3 we prove our next result which arose from the following questions we learned from L. Szpiro and T. Tucker. To proceed we need some definitions.

\bd \label{domination}
Let $(X, \phi), (Y,\psi)$ be dynamical systems polarized with respect to $\mc{L}, \mc{M}$ respectively.
We say $(X,\phi)$ \ti{dominates} $(Y,\psi)$ if there exists a dominant morphism $g: X \to Y$ such that:
\ben
\item $g \circ \phi = \psi \circ g$, and
\item $g^{*}\mc{M}= \mc{L}$
\een
\ed

Szpiro and Tucker asked the following:
\ben
\item Let $X,Y$ be projective $K$-varieties and $g: X \to Y$ a $K$-morphism. If $X$ is isotrivial, does it imply
that $Y$ is isotrivial?
\item Let $(X,\phi), (Y,\psi)$ be dynamical systems where the projective $K$-varieties $X,Y$ are isotrivial
and $g: X\to Y$ a dominant (or surjective) $K$-morphism
such that $g\circ \phi = \psi \circ g$. If $\phi$ is isotrivial, does it imply that $\psi$ is isotrivial?
\een

In \cite{B} the first author answered the first question affirmatively whenever $g_{*}\mc{O}_X = \mc{O}_{Y}$.
Our strategy to answer the second question in certain cases is to show that at
each place of $K$ the reduction of $\psi$ has degree $d^{\dim Y}$. This forces
reduction of $\psi$ to be well defined at every point.
Let us denote reduction with tildes. If $\tilde{\psi}$ has
potential good reduction at every place of $K$ there is nothing to show. Otherwise, we first prove:

\begin{reptheorem}{morphism}
Let $(X,\phi), (Y,\psi)$ be dynamical systems over $K$ polarized with respect to $\mc{L}, \mc{M}$ 
respectively, where $X,Y$ are connected. Assume that $X,Y, \phi$ are isotrivial, $\mc{L}, \mc{M}$ 
are defined over $k$, $(X,\phi)$ dominates $(Y,\psi)$ and one of the following conditions holds:
\ben
\item $\dim X = \dim\tilde{g}(X)$ and $\tilde{g}$ is dominant.
\item For some $n$ there exists a real number $d>0$ such that the action of $(\phi^{n})^{*}$ on $\Pic X$ is multiplication by $d^{n}$.
\item $Y = \mb{P}^{1}_{K}$.
\een
Then the rational map $\tilde{\psi}\mid_{\tilde{g}(X)}: \tilde{g}(X) \to \tilde{g}(X)$ is a morphism.
\end{reptheorem}

Now we have the following corollaries:

\begin{repcor}{domred}
Let $(X,\phi), (Y,\psi)$ be dynamical systems as in the previous theorem.
Suppose one of the conditions of Theorem \ref{morphism} holds, $\tilde{g}$ is
dominant at every place of $K$, and $\Aut(Y, \mathcal{M})$ is reductive. Then
$\psi$ is defined over some finite extension of $k$.
\end{repcor}

\bd
We say a non-empty proper subset $E\subset X$ is an \ti{exceptional set for $\phi$}
if $\phi^{-1}(E) =E=\phi(E)$.
\ed

\begin{repcor}{noexcep}
Let $(X,\phi), (Y,\psi)$ be dynamical systems as in the previous theorem.
Suppose one of the conditions of Theorem ~\ref{morphism} holds, $\Aut(Y,\mc{M})$
is reductive and $\phi$ has no non-empty proper exceptional subvarieties. Then
$\psi$ is defined over  some finite extension of $k$.
\end{repcor}

In the case $Y= \mb{P}^{1}_{K}$ we prove a stronger result using geometric invariant theory.
More precisely we prove,

\begin{repcor}{P1}
Let $(X,\phi), (\mb{P}^{1}_{K}, \psi)$ be dynamical systems polarized w.r.t. $\mc{L}, \mc{O}(1)$ 
respectively. Assume that $X$ is connected and $(X,\phi)$ dominates 
$(\mb{P}^{1}, \psi)$. Then $\psi$ is defined over some finite extension of $k$.
\end{repcor}


\section{Moduli Space of Polarized Dynamical Systems}

Let $k$ be an infinite field and $\Hom_{n,d}$ be the space of endomorphisms of
$\mb{P}^n_k$ of degree $d^n$.
By \S3.2 of \cite{PST}, $\Hom_{n,d}$ is an affine open subvariety of
$\mb{P}(\Sym^d((k^{n+1})^{n+1}))$, where $\Sym^d((k^{n+1})^{n+1})$ is the space of
homogeneous maps $k^{n+1} \to k^{n+1}$ of degree $d$.

Let $X$ be a connected, projective variety, and let $(X,\phi)$ be dynamical system over $k$ polarized w.r.t. 
$\mc{L}$. By a result of Fakhruddin \cite{F} (Thm. 2.1), any dynamical system $(X,\phi)$ polarized w.r.t. 
$\mc{L}$ extends to a dynamical system on projective space $(\mb{P}^n_k, \psi)$ polarized w.r.t. 
$\mc{O}(1)$. More precisely, we can find some $s \geq 1$ such that there exists a morphism $\psi$ on 
$\mathbb{P}^{n}$ making the following diagram commute
\begin{equation}
\label{polar}
\xymatrix@R+2em@C+2em{
X \ar[r]^{\phi} \ar[d]_{i_{\mc{L}^{\otimes s}}} & X \ar[d]^{i_{\mc{L}^{\otimes s}}} \\
\mathbb{P}^{n} \ar[r]^{\psi} & \mathbb{P}^{n} }
\end{equation}

\noindent In general, $\psi$ will not be unique. So we define the space of endomorphisms of $X$ as the quotient:
$$\Hom_d(X, \mc{L}) := \{ \psi \in \Hom_{n,d}: \psi (X) =X\}/ \{\psi \sim \psi' \iff
\psi|_X = \psi'|_X \}$$
as the analog of $\Hom_{n,d}$ for $\mathbb{P}^{n}$.

\begin{rem}
The construction of $\Hom_d(X,\mc{L})$ is similar to the construction of the space of polynomials in 
one dimension: the automorphism group we need to mod out by is the upper triangular group 
$\Gamma$, but we can instead embed the space of polynomials inside the space of rational maps, 
denoted by $\Rat_d$ and then study it as a subspace of the quotient $\Rat_d/\Gamma$. Of course we
cannot always mod out by a non-reductive group and expect a nice quotient, but the space of polynomials 
in one variable is sufficiently well-behaved: namely, we can consider \emph{centered} polynomials i.e. those 
whose $z^{d-1}$ coefficient is zero, and are acted on by a finite group of automorphisms. There is a natural 
bijection of sets between polynomials modulo the upper triangular group and centered polynomials modulo
the diagonal group.
\end{rem}

We also obtain (\emph{a priori} merely set-theoretic) quotients by the conjugation
action of automorphism group of $X$ associated to the line bundle class $\mc{L}$.
Observe that the full group $\Aut(X)$ does not act on $\Hom_d(X, \mc{L})$,
because an element of $\Aut(X)$ may map the line bundle class $\mc{L}$ to a different class; e.g. if $X$ is an elliptic curve, then the
translation group does not fix any ample class. Thus we only quotient by the subgroup of
automorphisms preserving the line bundle class i.e.
$$\Aut(X, \mc{L}) := \{\tau \in \Aut(X)\;|\; \tau^*\mc{L} \cong \mc{L} \}$$
To extend the GIT results in \cite{L1}, \cite{PST} to projective varieties, we need to show that
$$\M_d(X, \mc{L}) := \Hom_d(X, \mc{L})/\Aut(X,\mc{L})$$
is an affine geometric quotient (in the sense of \cite{M}). We prove this when
$\Aut(X,\mc{L})$ is reductive. To see that we can use GIT tools, we first need to show,

\bl
The space of endomorphisms of a polarized dynamical system, denoted by
$\Hom_{d}(X,\mc{L})$ is an affine scheme.
\el

\begin{proof}
The space $\Hom_{d}(X,\mc{L})$ arises as a closed subvariety of the space of morphisms
from $X$ to $\mb{P}^n$ induced by $\mc{L}^{\otimes s}$ for some integer $s >0$. Let $X'$
be the Zariski closure of the pullback of $X$ to $\mb{A}^{n+1} \setminus \{0\}$ i.e.
$$X' = \overline{ X \times_{\mb{P}^{n}} (\mb{A}^{n+1} \setminus \{0\})  } $$
A map fails to be a morphism if it is ill-defined at at least one point. Consider the space
$\Gamma(X,\mc{L}^{\otimes s}) \times X'$, which parametrizes homogeneous maps from
$X'$ to $\mb{A}^{n+1} \setminus\{0\}$ along with a point. Let $V$ be the subvariety defined by
$$ V = \{ (\psi,x) : \psi(x) =0 \} $$
i.e. the parameter space of maps that are ill-defined at some point. Now $V$ is closed. Moreover, for each $x \in X'$ there are $n+1$
sections that we equate to 0. Within the space of rational maps from $X$ to $X$ the
sections define a point on $\Gamma(X,\mc{L}^{\otimes s})$, so the condition of those sections
being all zero is actually of pure codimension at most $\dim X' = \dim X +1$.

Let $\pi_{1}$ be the projection from $\Gamma(X,\mc{L}^{\otimes s}) \times X'$
on the first factor.  The projection of $V$ onto the first factor has relative dimension $1$, 
since if $\psi$ is homogeneous and $\psi(x) = 0$ then $\psi(cx) = 0$ for all scalars $c$; 
its relative dimension is not higher because maps that are ill-defined at one point of $X$ 
are generically ill-defined at just one point of $X$. This means that 
$\dim\pi_{1}(V) = \dim V - 1$. Since $\codim V \leq \dim X'$, we have 
$\dim V \geq \dim\Gamma(X, \mc{L}^{\otimes s})$ and then the codimension of $\pi_{1}(V)$ 
is purely at most $1$. If there exist morphisms from $X$ to itself, then the codimension is actually 
$1$ and not $0$ and is again pure. Thus the complement is an affine open subvariety.
\end{proof}

Recall that the stabilizer of a point $\phi \in \Hom_d(X,\mc{L})$ under the action of
$\Aut(X,\mc{L})$ is the subgroup $\Stab(\phi) = \{ \gamma \in \Aut(X,\mc{L}) \; |\;
\gamma^{-1} \phi \gamma = \phi \}$. Note that if $P \in X$ is $\phi$-preperiodic for
some integers $l>m\geq1$, then $\gamma(P)$ is $\phi$-preperiodic for $l,m$ as well.
We use this fact in the next lemma.

\bl\label{finite}
For any $\phi \in \Hom_d(X,\mc{L})$ the stabilizer group, $\Stab(\phi)$ is finite
under the action of $\Aut(X,\mc{L})$.
\el

\bpf We follow the argument in \cite{PST} (Prop. 8). By \cite{F} (Thm. 5.1), $X$
has a dense set of preperiodic points for each $\phi$. In particular, there exists a
set $S$ of $n+2$ preperiodic points that span $\mathbb{P}^{n}$ having preperiod at
most $(l, m)$. Each $f \in \Stab(\phi)$ will act on the finite set of
preperiod-$(l, m)$ points, and moreover if $f$ acts trivially then it fixes the
spanning set $S$ and is therefore the identity. Thus $\Stab(\phi)$ embeds into a
finite group. \epf

\bl The action of $\Aut(X,\mc{L})$ on $\Hom_d(X,\mc{L})$ is closed. \el

\bpf
Following the argument in \cite{M} (p.10), if for each $\phi \in \Hom_d(X,\mc{L})$
there exists an open neighborhood $U$ of $\phi$ where the dimension of the
stabilizer is constant for all $\psi \in U$, then the action of $\Aut(X,\mc{L})$ on
$\Hom_d(X, \mc{L})$ is closed. Since $\Stab(\phi)$ is finite for all
$\phi \in \Hom_d(X,\mc{L})$, the action of $\Aut(X,\mc{L})$ on $\Hom_d(X,\mc{L})$ is closed.
\epf

\bt 
Let $X$ be a connected projective variety over an infinite field $k$ and let $\mc{L}$ be an ample 
line bundle on $X$. Then the scheme $\M_d(X,\mc{L})=\Hom_d(X,\mc{L})/\Aut(X,\mc{L})$ is an 
affine geometric quotient whenever $\Aut(X,\mc{L})$ is reductive.
\et

\bpf
By the previous lemma the action of $\Aut(X,\mc{L})$ on $\Hom_d(X,\mc{L})$ is
closed and $\Aut(X,\mc{L})$ is reductive by assumption. By Amplification 1.3 of
\cite{M} (p.30), $\Hom_d(X, \mc{L})/\Aut(X,\mc{L})$ is an affine geometric quotient.
Moreover it is equal to the spectrum of
$\Gamma(\Hom_d(X,\mc{L}), \mc{O}_{\Hom_d(X,\mc{L})})^{\Aut(X,\mc{L})}$.
\epf

As a corollary we obtain the following:

\bc\label{isotriv}
Let $K=k(C)$ be the function field of a smooth projective curve $C$ defined over
an infinite field $k$, $X$ a connected projective $k$-variety, $\mc{L}$ an ample line bundle on $X$,
$\phi: X_K \to X_K$ a polarized $K$-morphism. Assume that
the subgroup of automorphisms preserving the line bundle class
$$\Aut(X, \mc{L}):= \{ \tau \in \Aut(X)\; |\; \tau^*\mc{L} \cong \mc{L} \} $$
is reductive. Then $\phi$ is isotrivial if and only if $\phi$ has potential good
reduction at all places $v$ of $K$. \ec

\bpf The only if direction is clear. If $\phi$ has potential good reduction at
every place, then each valuation of the function field corresponds to a point in
the moduli space $\M_d(X,\mc{L})$. As in \cite{PST} we obtain a morphism from a
complete curve to an affine variety $\M_d(X,\mc{L})$, hence the image
is a point. The fiber of this point in $\M_d(X,\mc{L})$ contains a unique
$\Aut(X, \mc{L})(k)$-conjugacy class. It follows that $\phi$ coincides with
the base extension $\psi_K: X_K \to X_K$ for some $\psi$ in this
class i.e. $\phi$ is isotrivial.
\epf

We now give an example which shows that the condition $\Aut(X,\mc{L})$ being reductive is non-trivial.

\begin{example}\label{blowup}Let $X$ be the blowup of $\mb{P}^{2}$ at a single point
$x$. Note that $X$ is Fano (i.e. $-K_X$ is ample) and every automorphism of $X$ fixes an ample class, namely the anticanonical class. Moreover there is a unique
exceptional curve on $X$ and every automorphism of $X$ fixes it, thus every automorphism of $X$ arises from an automorphism of $\mb{P}^2$ fixing $x$. Clearly the
group of automorphisms of $\mb{P}^2$ fixing $x$ is not reductive since the unipotent radical is non-trivial. In particular, the unipotent radical is of the form
$$
\begin{pmatrix}
1 & 0 & 0 \\
0 & 1 & 0 \\
* & * & 1
\end{pmatrix}
$$
\end{example}

\begin{rem}
Let $X$ be as above. We show that polarized endomorphisms of $X$ exist. For example, any morphism of $\mb{P}^{2}$ which leaves
$x$ totally invariant. In fact, these are the only polarized morphisms of $X$: since the only pair of
irreducible curves with negative intersection is the exceptional curve intersected
with itself, the preimage of the exceptional curve must contain the exceptional
curve, i.e. the forward image of the exceptional curve is itself. Now take a curve
that does not intersect the exceptional curve; its preimage will also not intersect
the exceptional curve, but will intersect any other curve on $X$, and so the
preimage of the exceptional curve cannot contain anything except the exceptional
curve. This means we can descend every morphism on $X$ to a morphism on
$\mb{P}^{2}$, and then we can construct the space of morphisms on $X$ as the
subvariety of morphisms on $\mb{P}^{2}$ leaving a point totally
invariant.
\end{rem}

In fact, the argument in Example~\ref{blowup} can be extended to many other varieties, 
as long as we allow iterates of $\phi$. 

\bp \label{polyh} 
Let $X$ be a connected, smooth projective variety, let $U$ be the ample cone in $\NS(X) \otimes \mb{R}$. If 
$\overline{U}$ is polyhedral, then every polarized morphism has an iterate that acts on $\NS(X)$ 
as scalar multiplication; the minimal iterate we need to take depends only on $U$.
\ep

\bpf Let $\phi$ be a polarized morphism on $X$, and consider the action of $\phi^{*}$ on 
$\NS(X) \otimes \mb{R}$. The action is a linear map preserving $U$, hence
also $\overline{U}$. If $\overline{U}$ is polyhedral then $\phi^{*}$ acts on the faces of 
$\overline{U}$; since there are finitely many faces, some iterate $(\phi^{n})^{*}$ maps each such 
face to itself. Now we can intersect faces to obtain lines, which will also be preserved. Since the 
span of $U$ is all of $\NS(X) \otimes \mb{R}$, and $U$ is open, there exists a basis for 
$\NS(X) \otimes \mb{R}$ (since $\NS(X) \otimes \mb{R}$ is a finite dimensional vector space) with
the property that each basis element lies on a ray in $\overline{U}$ that is an eigenvector
for $(\phi^{n})^{*}$, and also $u \in U$ only if $u$ is a strictly positive linear combination of those 
elements. But now there exists a $u \in U$ that is an eigenvector for $(\phi^{n})^{*}$, whence all 
of the basis eigenvectors have the same eigenvalue and the action is scalar.\epf

All Fano varieties have polyhedral ample cones \cite{CPS}, as do Mori dream spaces
\cite{5A}, such as blowups of $\mb{P}^{n}$ at sufficiently general points. There exist Fano varieties 
that have non-reductive automorphism groups, even without being blowups. For examples we refer 
the reader to \cite{N} (Thm. 1.4, 1.5). On a more positive side, an arbitrary variety will not have any 
polarized morphisms. In particular, a compact complex manifold of general type (i.e. with maximal 
Kodaira dimension) does not admit any endomorphisms of degree greater than 1; see \cite{I} (Prop. 
10.10). As a result of the difficulty of finding endomorphisms of degree greater than $1$ on arbitrary 
projective varieties, we do not have a counterexample of a variety $X$ with a non-reductive 
automorphism group on which Lemma~\ref{finite} or Corollary~\ref{isotriv} is false.

In particular, if $X$ is a blowup of $\mb{P}^{n}$ with polyhedral ample cone, then after choosing 
a suitable $m$, for every $\phi: X \to X$ the map $\phi^{m}$ fixes each exceptional divisor, and 
we can descend it to a morphism on $\mb{P}^{n}$ for which all of the blown up points are totally 
invariant.

\rem The ample cone is not always polyhedral. For example, if $E$ is an elliptic curve without 
complex multiplication, then the ample cone of $E \times E$ has a quadratic condition \cite{BS}.


\section{Descent of Morphisms on Projective Varieties}

Throughout this section we have the following setup: Let $C$ be a smooth projective curve defined over an infinite
field $k$ and let $K=k(C)$ be the function field of  $C$. Let $(X,\phi), (Y,\psi)$ be dynamical systems over $K$
polarized with respect to $\mc{L}, \mc{M}$ respectively, with $X,Y$ connected. The varieties $X,Y$ and the morphism $\phi$ are
isotrivial, $\mc{L}, \mc{M}$ are defined over $k$, and $(X,\phi)$ dominates $(Y,\psi)$. We denote Zariski closure
of a set $S$ by $\overline{S}$. We fix a place $v\in K$ throughout and denote the reduction of a morphism $\rho$
at $v$ by $\tilde{\rho}$.

\bt \label{morphism}
Let $(X, \phi), (Y, \psi)$ be dynamical systems polarized with respect to $\mc{L}, \mc{M}$ respectively, where
$X, Y$ are connected.
Assume that $X, Y, \phi$ are isotrivial, $\mc{L}, \mc{M}$ are defined over $k$, $(X, \phi)$
dominates $(Y, \psi)$ and one of the following conditions holds:
\begin{enumerate}
\item $\dim X = \dim\tilde{g}(X)$ and $\tilde{g}$ is dominant.
\item For some $n$ there exists a $d>0$ such that the action of $(\phi^{n})^{*}$ on $\Pic X$ is
multiplication by $d^{n}$.
\item $Y = \mb{P}^{1}$.
\end{enumerate}
Then the rational map $\tilde{\psi}\mid_{\tilde{g}(X)}: \tilde{g}(X) \to \tilde{g}(X)$ is a morphism.
\et

\bpf Write $Z= \overline{\tilde{g}(X)}$. By restricting $\tilde{\psi}$ to $\tilde{g}(X)$, the morphism
$\tilde{\psi}$ is polarized in the
sense that the divisor class of the intersection $M \cdot Z$ pulls back to a $d^{th}$ power of itself. It suffices
to prove that the map has degree $d^{\dim Z}$ on $Z$. Then it is a morphism on $Z$, in which case it is a
morphism on its open (not necessarily proper) subvariety $\tilde{g}(X)$ also.

The degree of $\tilde{\psi}\mid_{Z}$ is at most $d^{\dim Z}$. Let $z \in \tilde{g}(X)$, and consider the fiber
$\tilde{g}_{z}:= \tilde{g}^{-1}(z)$. The preimage of $\tilde{g}_{z}$ under $\phi$ consists of some fibers
$\tilde{g}_{z_{i}}$. There is a generic cycle class for fibers of $\tilde{g}$, of codimension $\dim\tilde{g}(X)$.
If we can show that this cycle class $C$ satisfies the equation $\phi^{*}(C) = d^{\codim C}C$ then we are
done, since we could pick $z$ to be suitably generic so that $\tilde{g}_{z}$ would have exactly
$d^{\dim\tilde{g}(X)}$ pre-image fibers.

If the first condition in the lemma holds, then $\dim C = 0$ and we know that each point has $d^{\dim X}$ preimages, counted with multiplicity (which is
generically $1$); see \cite{PST} (Prop. 2). If the second condition holds for $n = 1$, then $\phi$ acts on each codimension-$i$ cycle group by scalar multiplication by
$d^{i}$ and we are also done; if it holds for some $n$, then clearly $\tilde{\psi^{n}}\mid_{Z}$ has degree $d^{n\dim\tilde{g}(X)}$ and then
$\tilde{\psi}\mid_{Z}$ has degree $d^{\dim\tilde{g}(X)}$. Finally, suppose that the third condition holds. We can set $\mathcal{M} = \mathcal{O}(1)$, and then
$\tilde{g}_{z}$ is a section of $\mathcal{L}$ since we can enlarge $\tilde{g}$ so that its indeterminacy locus has codimension at least $2$, and $\phi^{*}$ acts on
$\mathcal{L}$ as scalar multiplication by $d$ by definition.\epf

\begin{rem} 
The second condition in Theorem~\ref{morphism} is automatically satisfied when $X$ has polyhedral ample cone, e.g. whenever $X$ is Fano.
\end{rem}

We now have the following corollaries:

\bc \label{domred}
Let $(X,\phi), (Y,\psi)$ be as above.
Suppose one of the conditions of Theorem \ref{morphism} holds, $\tilde{g}$ is
dominant at every place of $K$, and $\Aut(Y, \mathcal{M})$ is reductive.
Then $\psi$ is defined over some finite extension of $k$.
\ec

\bpf If $\tilde{g}$ is dominant, then $\dim\tilde{g}(X) = \dim Y$, and then
$\tilde{\psi}$ has degree $d^{\dim Y}$, making it a morphism on $Y$. In other
words, $\psi$ has good reduction at every place. Then by Corollary~\ref{isotriv},
$\psi$ is defined over some finite extension of $k$.\epf

\begin{rem} 
In the previous corollary we require $k$ to be infinite since we invoke Corollary \ref{isotriv}.
\end{rem}

\bl \label{excep}
Let $(X,\phi), (Y,\psi)$ be as above. Suppose one of the conditions
of Theorem ~\ref{morphism} holds. Then the locus of indeterminacy of $\tilde{g}$
is an exceptional set for $\phi$, i.e. it is its own preimage.
\el

\bpf Let $W$ be the indeterminacy locus of $\tilde{g}$ and $Z$ be the indeterminacy locus of $\tilde{\psi}\mid_{\tilde{g}(X)}$. The indeterminacy locus of
$\tilde{g}\circ\phi$ is then $\phi^{-1}(W)$ and the indeterminacy locus of $\tilde{\psi}\circ{\tilde{g}}$ is $W \cup \tilde{g}^{-1}(Z)$. Since
$\tilde{\psi}\circ\tilde{g} = \tilde{g}\circ\phi$, we obtain that $W \cup \tilde{g}^{-1}(Z) = \phi^{-1}(W)$. Now we apply Theorem \ref{morphism} and obtain that $Z$
is empty, so that $W = \phi^{-1}(W)$.\epf

In general, there are no exceptional subvarieties for a morphism on $\mb{P}^{n}$. Conjecturally (\cite{FN}, Conjecture 2.6) all such subvarieties are linear, and in the linear case, a large
number of monomials of $\phi = (\phi_{0}:\ldots:\phi_{n})$ has to be zero for $\phi$ to have an exceptional subvariety. This motivates our following corollary:


\bc \label{noexcep}
Let $(X,\phi), (Y,\psi)$ be as above.
Suppose one of the conditions of Theorem ~\ref{morphism} holds, $\Aut(Y,\mc{M})$ is reductive
and $\phi$ has no non-empty proper exceptional subvarieties. Then $\psi$ is defined over some 
finite extension of $k$.
\ec

\bpf By Corollary~\ref{excep}, $\tilde{g}$ is a morphism from $X$ to $Y$. If we can show that $\tilde{g}$ is
dominant, then by Corollary~\ref{domred} we are done.

If $X = Y$ and $g$ is polarized then by Corollary~\ref{isotriv}, $g$ is isotrivial, and $\tilde{g} = g$ is dominant (in fact surjective).
Otherwise, we do not have an equivalent of Corollary~\ref{isotriv}, but we do have something close enough. A rational map from $X$ to $Y$ may be viewed as a
subvariety of $X \times Y$ (i.e. the graph of $g$, denoted $\Gamma_{g}$) that intersects the cycle class $\{x\} \times Y$ in exactly one point with multiplicity $1$; it is a morphism if and only if
$\Gamma_{\tilde{g}}$ intersects each subvariety $\{x\} \times Y$ in one point with multiplicity $1$, rather than just the generic cycle class. Now consider the intersection $\Gamma_{g} \cdot (X \times \{y\})$. Since $g$ is surjective, the intersection $\Gamma_{g} \cdot (X \times \{y\})$ is a nonzero cycle class of dimension $\dim X - \dim Y$. Finally, the intersection
$\Gamma_{\tilde{g}} \cdot (X \times \{y \})$ is still a nonzero cycle class; this means that for the generic $y \in Y$, $(X \times \{y\}) \cdot \Gamma_{\tilde{g}}$ is nonempty, and since $\tilde{g}$
is a morphism, we have $(x, y) \in \Gamma_{\tilde{g}}$ if and only if $y = \tilde{g}(x) \in \tilde{g}(X)$. Thus we get that $\overline{\tilde{g}(X)} = Y$ (in fact,
$\tilde{g}(X) = Y$) and $\tilde{g}$ is dominant (in fact, surjective).\epf

In the case $Y= \mb{P}^{1}_{K}$, we give another proof to show that $\psi$ is defined over 
some finite extension $k$. We use the notion of semistability from geometric invariant theory 
to prove this.

\bd \label{hmc}(Hilbert-Mumford Criterion) Let $G$ be a geometrically reductive group acting 
on a projective variety $X$, with $\mathcal{L}$ a $G$-invariant line  bundle on $X$. By embedding 
$X$ into $\mb{P}^{n}$ via $\mathcal{L}$, we obtain an action of $G$ on $\mb{P}^{n}$, so we may 
assume that $X = \mb{P}^{n}$. For each one-parameter subgroup $\mathbb{G}_{m} \subseteq G$, 
pick a coordinate system for $\mb{P}^{n}$ with respect to which its action is diagonal, and consider 
the weight $t_{i}$ with which the subgroup acts on each coordinate $x_{i}$. A point $x \in X$ is said 
to be \ti{unstable} (respectively \ti{not stable}) if there exists a one-parameter subgroup such that, 
after the appropriate coordinate change, we have $x_{i} = 0$ whenever $t_{i} \leq 0$ 
(resp. $t_{i} < 0$); otherwise, it is said to be \ti{semistable} (resp. \ti{stable}).\ed

Recall from geometric invariant theory the following theorem:

\bt \label{ssrorig}Let $G$ be a geometrically reductive group acting on a projective variety $X$ whose stable and semistable points are denoted by  $X^{s}$ and
$X^{ss}$ respectively. Let $R$ be a discrete valuation ring with fraction field $K$, and let $x_{K} \in X^{s}_{K}$. Then for some finite extension $K'$ of $K$, with
$R'$ the integral closure of $R$ in $K'$, $x_{K}$ has an integral model over $R'$ with semistable reduction modulo the maximal ideal. In other words, we can find
some $A \in G(\overline{K})$ such that $A\cdot x_{K}$ has semistable reduction. If $x_{K} \in X^{ss}_{K}$, then the same result is true, except that $x_{R'}$ could
be an integral model for some $x'_{K'}$ mapping to the same point of $X^{ss}//G$ such that $x'_{K'} \notin G\cdot x_{K}$.\et

\bpf \cite{L2} (Theorem 2.11). \epf

\bd
Let $K$ be a complete non-archimedean field with ring of integers $R$, and let 
$\Hom(X_K,Y_K)$ be the space of $K$-morphisms from $X_K$ to $Y_K$. We say 
$\phi_R \in \Hom(X_R, Y_R)$ is an \ti{integral model} of $\phi_K \in  \Hom(X_K,Y_K)$ if 
$\phi_R$ is generically equal to $\phi_K$.
\ed

\bl \label{ssr} Let $X$ be a projective variety over $k$ and $g: X_K \to \mb{P}^{n}_{K}$ be a $K$-morphism, where $K$ is a complete non-archimedean field. Then $g$
has an integral model whose reduction is non-constant. \el

\bpf By Theorem~\ref{ssrorig}, it suffices to show that constant maps are unstable with respect to some action.
Even though $\Aut\mb{P}^{n} = \PGL_{n+1}$;
for the purposes of the Hilbert-Mumford criterion, $\SL_{n+1}$ is the correct group to use. Note that we are
allowed to take finite extensions, and $\SL_{n+1}$ maps finite-to-one onto $\PGL_{n+1}$. We choose the left action
of $\SL_{n+1}$ on the space of maps from $X$ to $\mb{P}^{n}$ i.e. for $h \in \SL_{n+1}, g: X \to \mb{P}^{n}$, we
have $h\cdot g(x) = h\circ g$.

Now, let $T$ be a one-parameter subgroup of $\SL_{n+1}$, say with diagonal entries $t^{a_{0}}, \ldots, t^{a_{n}}$,
with $$a_{0} \geq \ldots \geq a_{n} \quad \tr{and} \quad \sum a_{i} = 0$$

A morphism from $X$ to $\mb{P}^{n}$ is given by an $(n+1)$-tuple of sections of some line bundle on $X$, say $g_{0}, \ldots, g_{n}$, and the action of $T$ has
weight $a_{i}$ on each $g_{i}$. Now if the morphism is constant, or even has image contained in a hyperplane, we assume the image is contained in $x_{n} = 0$, so
that $g_{n} = 0$. If $a_{0} = \ldots = a_{n-1} = 1$ and $a_{n} = -n$ then the coordinates of $g$ are zero whenever the action of $T$ has negative weights
and thus $g$ is unstable.

Finally, observe that, with respect to the $\mathcal{O}(1)$-bundle on $\mb{P}^{n}$, maps whose images are not contained in a hyperplane are not just semistable but
also stable. This is because $T$ acts on each $g_{i}$ with weight depending only on $i$; therefore, if for some weight all coordinates of $g$ are zero then we
have $g_{i} = 0$ and then $g$ is unstable. Thus the only group that can act with no negative weights but with some non-negative weights is the group all of whose
weights are zero, i.e. the trivial group. This means that $g$ is stable, and thus it has an integral model with semistable (in fact, stable) reduction.\epf

\bc \label{P1}
Let $(X,\phi), (\mb{P}^{1}_{K}, \psi)$ be dynamical systems polarized w.r.t. $\mc{L}, \mc{O}(1)$ 
respectively. Assume that $X$ is connected and $(X,\phi)$ dominates 
$(\mb{P}^{1}, \psi)$. Then $\psi$ is defined over some finite extension of $k$.
\ec

\bpf
A morphism from a connected variety $X$ to $\mb{P}^{1}_{K}$ is either constant or dominant. After conjugation we can force $\tilde{g}$ to be non-constant by Lemma~\ref{ssr}. Now apply Corollary~\ref{domred}.
\epf

\ti{Funding.} The first author was partially supported by PSC-CUNY grant TRADA-44-179. \\

\ti{Acknowledgements.} 
This research project was initiated while the authors were postdoctoral fellows at ICERM, Brown
University. We would like to thank Joseph Silverman and the staff at ICERM for their 
hospitality, and for providing an environment conducive for our research.

\end{document}